\documentclass{amsart}
\usepackage{amsthm}
\usepackage{enumerate}
\newtheorem{thm}{Theorem}
\newtheorem{prop}[thm]{Proposition}
\newtheorem{rem}[thm]{Remark}
\newtheorem{lem}[thm]{Lemma}
\newtheorem{cor}[thm]{Corollary}

\newtheorem{ex}[thm]{Example}
\newcommand{\conv}{\mathrm{Conv}}
\newcommand{\extr}{\mathrm{extr}}
\newcommand{\supp}{\mathrm{supp}}
\newcommand{\dist}{\mathrm{dist}}
\newcommand{\mmod}{\mathrm{mod}}

\title{Orbits in symmetric spaces}
\author{F. Sukochev and D. Zanin}

\begin{document}
\begin{abstract} We characterize those elements in a fully symmetric spaces on the interval $(0,1)$ or on the semi-axis $(0,\infty)$ whose orbits are the norm-closed convex hull of their extreme points.
\end{abstract}

\maketitle

\section{Introduction}

The following semigroups of bounded linear operators  play a fundamental role in the interpolation theory of linear operators for the couple $(L_1, L_\infty)$ of Lebesgue measurable functions on intervals $(0,1)$ and $(0,\infty)$. The semigroup of absolute contractions, or admissible operators (see e.g. \cite[II.3.4]{KPS})
$$\Sigma:=\{ T:L_1+L_{\infty}\to L_1+L_{\infty}:\ \max(\|T\|_{L_1\to L_1},\|T\|_{L_{\infty}\to L_{\infty}})\le 1,$$
the semigroup of substochastic operators (see e.g. \cite[p.107]{BS})
$$\Sigma_+:=\{0\leq T\in\Sigma\}$$  and,
in the case of the interval $(0,1)$, the semigroup of doubly stochastic operators
$$\Sigma':=\{0\leq T\in\Sigma_+:\ \int_0^1(Tx)(s)ds=\int_0^1x(s)ds,\ \forall  x\ge 0,\  T1=1\}$$
(see e.g. \cite{Ryfftrans}). If $x\in L_1+L_\infty$ (respectively, $0\leq x\in L_1+L_\infty$ or $0\leq x\in L_1(0,1)$) we denote by $\Omega(x)$ (respectively $\Omega_+(x)$ and $\Omega'(x)$) the {\it orbit} of $x$ with respect to the semigroups $\Sigma$ (respectively, $\Sigma_+,$ and $\Sigma'$). A Banach function space $E$ (on  $(0,1)$ or $(0,\infty)$, see \cite[pp.2-3]{BS}) is called an exact interpolation space if every $T\in \Sigma$ maps $E$ into itself and $\|T\|_{E\to E}\leq 1$, or alternatively, if $\Omega(x)\subset E$  and $\|y\|_E\leq \|x\|_E$ for every $x\in E$ and every $y\in\Omega(x).$ The class of exact interpolation spaces admits an equivalent description in terms of (sub)majorization in the sense of Hardy, Littlewood and Polya.  Recall, that if $x,y\in L_1+L_{\infty},$ then $y$ is said to be a {\it submajorized} by $x$ in the sense of Hardy, Littlewood and Polya, written $y\prec\prec x$  if and only if
$$\int_0^ty^*(s)ds\leq\int_0^tx^*(s)ds\quad  t\geq 0.$$
Here, $x^*$ denotes the non-increasing right-continuous rearrangement of $x$ given by
$$x^*(t)=\inf\{s\geq0:\ m(\{|x|\geq s\})\leq t\}$$ and $m$ is Lebesgue measure.
If $0\leq x,y\in L_1$, then we say that $y$ is {\it majorized} by $x$ (written $y\prec x$) if and only if $y\prec\prec x$ and $||y||_1=||x||_1$. A Banach function space $E$ is said to be {\it fully symmetric} if and only if $x\in E,$ $y\in L_1+L_{\infty}$ $y\prec\prec x$ $\Rightarrow$ $y\in E$ and $||y||_E\leq||x||_E$.
The classical Calderon-Mityagin theorem (see \cite{Calderon1966}, \cite{KPS}, \cite{BS}) gives an alternative description of the sets $\Omega(x)$, $x\in L_1+L_\infty$ and $\Omega_+(x)$, $0\leq x\in L_1+L_\infty$ as follows
$$\Omega(x)=\{y\in L_1+L_{\infty}:\ y\prec\prec x\},\ \Omega_{+}(x)=\{0\leq y\in L_1+L_{\infty}:\ y\prec\prec x\}$$
and (in the case of the interval $(0,1)$ and $0\leq x\in L_1(0,1)$)
$$\Omega'(x)=\{0\leq y\in L_1:\ y\prec x\},$$
which shows, in particular, that the classes of exact interpolation spaces and fully symmetric spaces coincide.

Let fully symmetric Banach function space $E$ be fixed. The principal aim of the paper is to give conditions for a given $0\leq x\in E$ which are necessary and sufficient for each of the sets $\Omega_+(x)$, $\Omega'(x)$ to be the norm closure of the convex hull of their extreme points.

If $E=L_1(0,1),$ then it has been shown by Ryff (see \cite{Ryfftrans}) that if $0\leq x\in E,$ then the orbit $\Omega'(x)$ is weakly compact and so, due to the Krein-Milman theorem, the orbit $\Omega'(x)$ is the weak (and hence norm)-closed convex hull of its extreme points. It follows from the results of \cite{DSS} that the set $\Omega'(x)$ is weakly compact in any separable symmetric space $E.$ Hence, $\Omega'(x)$ is the weak (and hence norm)-closed convex hull of its extreme points in any separable symmetric space $E.$

If a fully symmetric space $E$ is not separable, then it is not the case in general that orbits are weakly compact. A trivial example yields the orbit $\Omega(\chi_{[0,1]})$ in fully symmetric non-separable space $L_\infty(0,1)$. Indeed, it is obvious that this orbit coincides with  the unit ball of $L_\infty(0,1)$ and the latter is not weakly compact since the space $L_\infty(0,1)$ is non-reflexive. Nonetheless, it is an interesting question to give necessary and sufficient conditions that the orbit of a given element should be the norm-closed convex hull of its extreme points. This question was considered by Braverman and Mekler (see \cite{BM}) in the case of the interval $(0,1)$ and orbits $\Omega(x)$. They showed that for every fully symmetric space $E$ on  $(0,1)$ satisfying the condition
\begin{equation}\label{bm condition}
\lim_{\tau\to\infty}\frac1{\tau}||\sigma_{\tau}||_{E\to E}=0
\end{equation}
that $\Omega(x)$ is indeed the norm-closed convex hull of the set of its extreme points, for every $x\in E$ (see \cite[Theorem 3.1]{BM}). Here $\sigma_{\tau}$ denotes the usual dilation operator (see the following section for definition and properties). They showed as well that the converse assertion is valid in case that $E$ is a Marcinkiewicz space on $(0,1)$. As explained below, this converse assertion, however, fails for arbitrary fully symmetric spaces.

We show (Theorem \ref{qprime finite}) that if $E$ is a fully symmetric space on $(0,1)$ and if $0\leq x\in E,$ then $\Omega'(x)$ is the norm-closed convex hull of its extreme points if and only if
\begin{equation}\label{new condition}
\varphi(x):=\lim_{\tau\to\infty}\frac1{\tau}||\sigma_{\tau}(x^*)||_E=0.
\end{equation}
As shown in Corollary \ref{result of BM} this implies the result of Braverman and Mekler. In the Appendix, we demonstrate that the conditions \eqref{bm condition} and \eqref{new condition} are distinct in the class of Orlicz spaces. If $E$ is an Orlicz space, then it is the case that \eqref{new condition} holds, and so by Theorems \ref{qprime finite} and \ref{qplus finite} for every $0\leq x\in E,$ the sets $\Omega'(x)$, $\Omega_+(x)$ and $\Omega (x)$ are the norm-closed convex hulls of its extreme points. However, there are non-separable Orlicz spaces $E$ which fail condition \eqref{bm condition}.

In the Appendix, we also introduce the notion of symmetric and fully symmetric functionals. The latter are a ``commutative'' counterpart of Dixmier traces appeared in non-commutative geometry (see e.g. \cite{CaSu}). Symmetric and fully symmetric functionals are extensively studied in \cite{DPSS1998}, \cite{KS} (see also \cite{CaSu} and references therein). Note, however, that our terminology differs from that used in just cited articles. A subclass of Marcinkiewicz spaces admitting symmetric functionals which fail to be fully symmetric is described in \cite{KS}. It follows from our results that any symmetric functional on a fully symmetric space satisfying \eqref{new condition} is automatically fully symmetric. In particular, this implies that an Orlicz space does not possess any singular symmetric functionals (see Proposition \ref{result on symm_funct}). This latter result strengthens the result of \cite[Theorem 3.1]{DPSS1998} that an Orlicz space does not possess any singular fully symmetric functionals.

Results similar to Theorems \ref{qprime finite} and \ref{qplus finite} hold also for fully symmetric spaces $E$ on the semi-axis (see Theorems \ref{qprime integrable}, \ref{qplus integrable}, \ref{qprime nonintegrable}, \ref{qplus nonintegrable}).

The main results of this article are contained in Section 4. In the following section we present some definitions from the theory of symmetric spaces, as some of our results hold in a slightly more general setting than that of fully symmetric spaces. For more details on the latter theory we refer to \cite{KPS, LT2, BS}. Section 3 treats various properties of the functional $\varphi$ and its modifications needed in Section 4. We would like to emphasize the difference between geometric properties of the orbit $\Omega(x)$ and those of $\Omega'(x)$ and $\Omega_+(x)$. This is especially noticeable in the description of the respective sets of their extreme points. The extreme points of the sets $\Omega(x)$ and $\Omega'(x)$, $x\ge 0$ are well-known (see \cite{Ryffproc,CKS}) and are given by:
$$\extr(\Omega(x))=\{y\in L_1+L_{\infty}:\ y^*=x^*\},\quad \extr(\Omega'(x))=\{0\leq y\in L_1:\ y^*=x^*\}.$$
whereas the description of extreme points of the set $\Omega_+(x)$, $x\ge 0$ given by
$$\extr(\Omega_{+}(x))=\{0\leq y\in L_1+L_{\infty}:\ y^*=x^*\chi_{[0,\beta]}\ for\ some\ \beta\leq\infty\}$$
when $x^*(\infty):=\lim _{t\to \infty}x^*(t)=0,$ and by
$$\extr(\Omega_{+}(x))= \{0\leq y\in L_1+L_{\infty}:\
y^*=x^*\chi_{[0,\beta]}\ for\ some\ \beta\leq \infty$$
$$and \ y\chi_{\{y<y^*(\infty)\}}=0\}$$
when $x^*(\infty)>0,$ is somewhat less known, so we present in the Appendix a careful exposition of the latter equalities.

{\bf Acknowledgements.}  We would like to thank Peter Dodds for many helpful comments on the content of this paper and lengthy discussions of earlier drafts. We also thank Sergei Astashkin for his interest.

\section{Preliminaries}

Let $L_0$ be a space of Lebesgue measurable functions either on $(0,1)$ or on $(0,\infty)$ finite almost everywhere (with identification $m-$a.e.). Here $m$ is Lebesgue measure. Define $S_0$ as the subset of $L_0$ which consists of all functions $x$ such that $m(\{|x|>s\})$ is finite for some $s.$

Let $E$  be a Banach space of real-valued Lebesgue measurable functions either on $(0,1)$ or $(0,\infty)$ (with identification $m-$a.e.). $E$ is said to be {\it ideal lattice} if $x\in E$ and $|y|\leq |x|$ implies that $y\in E$ and $||y||_E\leq||x||_E.$

The ideal lattice $E\subseteq S_0$ is said to be {\it symmetric space} if for every $x\in E$ and every $y$ the assumption $y^*=x^*$ implies that $y\in E$ and $||y||_E=||x||_E.$


If $E=E(0,1)$ is a symmetric space on $(0,1),$ then
$$L_{\infty}\subseteq E\subseteq L_1.$$

If $E=E(0,\infty)$ is a symmetric space on $(0,\infty),$ then
$$L_1\cap L_{\infty}\subseteq E\subseteq L_1+L_{\infty}.$$


Symmetric space $E$ is said to be {\it fully symmetric} if and only if $x\in E,$ $y\in L_1+L_{\infty}$ $y\prec\prec x$ $\Rightarrow$ $y\in E$ and $||y||_E\leq||x||_E.$

We now gather some additional terminology from the theory of symmetric spaces that will be needed in the sequel.

Suppose $E$ is a symmetric space. Following \cite{BM}, $E$ will be called {\it strictly symmetric} if and only if whenever $x,y\in E$ and $y\prec\prec x$ then $||y||_E\leq||x||_E.$

It is clear that if $E$ is fully symmetric then $E$ is strictly symmetric, but the converse assertion is not valid.

The norm $||\cdot||_E$ is called Fatou norm if, for every sequence $x_n\uparrow x\in E,$ it follows that $||x_n||_E\uparrow ||x||_E$. This is equivalent to the assertion that the unit ball of $E$ is closed with respect to almost everywhere convergence.

It is well known that if the norm on $E$ is a Fatou norm then $E$ is strictly symmetric.

If $\tau>0,$ the dilation operator $\sigma_{\tau}$ is defined by setting $(\sigma_{\tau}(x))(s)=x(\frac{s}{\tau}),$ $s>0$ in the case of the semi-axis. In the case of the interval $(0,1)$ the operator $\sigma_{\tau}$ is defined by
$$
(\sigma_{\tau}x)(s)=
\begin{cases}
x(s/\tau),& s\leq\min\{1,\tau\}\\
0,& \tau<s\leq1.
\end{cases}
$$

The operators $\sigma_{\tau}$ ($\tau\geq1$) satisfy semi-group property $\sigma_{\tau_1}\sigma_{\tau_2}=\sigma_{\tau_1\tau_2}.$ If $E$ is a symmetric space and if $\tau>0,$ then the dilation operator $\sigma_{\tau}$ is a bounded operator on $E$ and
$$||\sigma_{\tau}||_{E\to E}\leq\max\{1,\tau\}.$$

We will need also the notion of a partial averaging operator (see \cite{BM}).

Let $\mathcal{A}=\{A_k\}$ be a (finite or infinite) sequence of disjoint sets of finite measure and denote by $\mathfrak{A}$ the collection of all such sequences. Denote by $A_{\infty}$ the complement of $\cup_kA_k.$ The partial averaging operator is defined by
$$P(x|\mathcal{A})=\sum_k\frac1{m(A_k)}(\int_{A_k}x)\chi_{A_k}+x\chi_{A_{\infty}}.$$
Note, that we do not require $A_{\infty}$ to have a finite measure.

Every partial averaging operator is a contraction both in $L_1$ and $L_{\infty}.$ Hence, $P(\cdot|\mathcal{A})$ is also contraction in $E.$ In case of the interval $(0,1),$ $P(\cdot|\mathcal{A})$ is a doubly stochastic operator in the sense of \cite{Ryfftrans}.

Since $P(\cdot|\mathcal{A})\in\Sigma,$ then $P(x|\mathcal{A})\in\Omega(x)$ (respectively, $P(x|\mathcal{A})\in\Omega'(x)$ if $x\in L_1$) for every $\mathcal{A}\in\mathfrak{A}.$ As will be seen, elements of the form $P(x|\mathcal{A})$ play a central role.

The following properties of rearrangements can be found in \cite{KPS}.
If $x,y\in L_1+L_{\infty},$ then
\begin{equation}\label{first property of rearrangement}
(x+y)^*\prec\prec x^*+y^*
\end{equation}
and
\begin{equation}\label{second property of rearrangement}
(x^*-y^*)\prec\prec(x-y)^*.
\end{equation}

Let us recall some classical examples of fully symmetric spaces.

Let $\psi$ be a concave increasing continuous function. The Marcinkiewicz space $M_{\psi}$ is the linear space of those functions $x\in S_0,$ for which
$$||x||_{M_{\psi}}=\sup_t\frac1{\psi(t)}\int_0^tx^*(s)ds<\infty$$
Equipped with the norm $||x||_{M_{\psi}},$ $M_{\psi}$ is a fully symmetric space with Fatou norm.

Let $M(t)$ be a convex function on $[0,\infty)$ such that $M(t)>0$
for all $t>0$ and such that
\begin{equation}\label{convex function}0=M(0)=\lim\limits_{t\to 0}\frac{M(t)}{t}=\lim\limits_{t\to
\infty}\frac{t}{M(t)}
\end{equation}
Denote by $L_M$ the Orlicz space on $[0,\infty)$ (see e.g.
~\cite{LT2,KPS}) endowed with the norm
\[\|x\|_{L_M}=\inf\{\lambda :\lambda >0,\;\int\limits_{0}^{\infty}
M(|x(t)|/\lambda)dt\leq 1\}.\]
Equipped with the norm $\|x\|_{L_M},$ $L_M$ is a fully symmetric space with Fatou norm.

For further properties of Marcinkiewicz and Orlicz spaces, we
refer to ~\cite{KPS, LT2} and ~\cite{RaRe1991}.

For $0\leq x\in L_1+L_{\infty},$ we set
$$Q_+(x)=\overline{\conv}(\extr(\Omega_+(x))).$$
For $0\leq x\in L_1,$ we set
$$Q'(x)=\overline{\conv}(\extr(\Omega'(x))).$$
For $0\leq x\in L_1+L_{\infty},$ we set
$$Q'(x)=\overline{\conv}\{y^*=x^*,\ y\chi_{\{y<y^*(\infty)\}}=0\}.$$
Here, $\overline{\conv}$ means the norm-closed convex hull. See
Appendix for the precise description of the extreme points.

\section{ The dilation functional and its properties}

The following assertion is widely used in literature. However, no direct reference is available.
\begin{lem}\label{well-known property} If $0\leq x,y\in L_1+L_{\infty},$ then  \begin{equation}\label{third property of rearrangement}
x^*+y^*\prec\prec2\sigma_{\frac12}((x+y)^*).
\end{equation}
\end{lem}

\begin{proof} Fix $\varepsilon>0.$ It follows from \cite[II.2.1]{KPS},
$$\int_0^tx^*(s)ds\leq\varepsilon+\int_{e_1}x(s)ds,\quad\int_0^ty^*(s)ds\leq\varepsilon+\int_{e_2}y(s)ds$$
for some $e_1$ and $e_2$ with $m(e_i)=t.$ However,
$$\int_{e_1}x(s)ds+\int_{e_2}y(s)ds\leq\int_{e_1\cup e_2}(x+y)(s)ds\leq$$
$$\leq\sup_{m(e)=2t}\int_e(x+y)(s)ds=\int_0^{2t}(x+y)^*(s)ds.$$
Note, that $\int_0^{2t}u(s)ds=\int_0^t(2\sigma_{\frac12}u)(s)ds.$
\end{proof}

\begin{lem}\label{dilation property} If $x,y\in L_1+L_{\infty}$ and $y\prec\prec x,$ then,
$$(\sigma_{\tau}(y))^*\leq\sigma_{\tau}(y^*)\prec\prec\sigma_{\tau}(x^*).$$
\end{lem}
\begin{proof} Set $d_y(s)=m(t:|y(t)|>s).$ In the case of the semi-axis, $d_{\sigma_{\tau}y}=\tau d_y=d_{\sigma_{\tau}(y^*)}.$ In the case of the interval $(0,1),$ $d_{\sigma_{\tau}y}\leq\tau d_y$ and $d_{\sigma_{\tau}(y^*)}=\min\{1,\tau d_y\}.$ Hence, $d_{\sigma_{\tau}y}\leq d_{\sigma_{\tau}(y^*)}$ and so $(\sigma_{\tau}(y))^*\leq\sigma_{\tau}(y^*).$ Finally,
$$\int_0^t\sigma_{\tau}(y^*)(s)ds=\tau\int_0^{\frac{t}{\tau}}y^*(s)ds\leq\tau\int_0^{\frac{t}{\tau}}x^*(s)ds=\int_0^t\sigma_{\tau}(x^*)(s)ds.$$
\end{proof}

The next lemma introduces the dilation functional $\varphi$ on $E,$ which is a priori non-linear. The behavior of the functional $\varphi$ on the positive part $E_+$ of $E$ provides the key to our main question.

\begin{lem}\label{general properties of fi} For every $x\in E$ the following limit exists and is finite.
\begin{equation}\label{definition of fi}
\varphi(x)=\lim_{s\rightarrow\infty}\frac1s||\sigma_s(x^*)||_E,\ x\in E.
\end{equation}
If, in addition, $E=E(0,\infty)$, then the following limits exist and are finite.
\begin{equation}\label{definition of fi fin}
\varphi_{fin}(x)=\lim_{s\rightarrow\infty}\frac1s||\sigma_s(x^*)\chi_{[0,1]}||_E,\ x\in E,
\end{equation}
\begin{equation}\label{definition of fi cut}
\varphi_{cut}(x)=\lim_{s\rightarrow\infty}\frac1s||\sigma_s(x^*)\chi_{[0,s]}||_E,\ x\in E.
\end{equation}
The following properties hold.
\begin{enumerate}[i)]
\item\label{zero property fi} If $E$ symmetric, then $\varphi(y)\leq\varphi(x)$ provided that $x,y\in E$ satisfy $ y^*\leq x^*$.
\item\label{first property fi}  If $E$ symmetric, then $\varphi(x)\leq||x||_E$ for every $x\in E.$
\item\label{second property fi} If $E$ strictly symmetric, then $\varphi(y)\leq\varphi(x)$ provided that $x,y\in E$ satisfy $y\prec\prec x.$
\item\label{third property fi} If $E$ is symmetric, then $\varphi(\sigma_{\tau}(x^*))=\tau\varphi(x),$ $\tau>0.$
\item\label{fourth property fi} If $E$ is strictly symmetric, then $\varphi$ is norm-continuous.
\item\label{fifth property fi} If $E$ is strictly symmetric, then $\varphi$ is convex.
\end{enumerate}
If, in addition, $E=E(0,\infty),$ then $\varphi_{fin}$ also satisfies \eqref{zero property fi},\eqref{first property fi},
\eqref{second property fi},\eqref{third property fi}, \eqref{fourth property fi} and \eqref{fifth property fi}, while $\varphi_{cut}$ satisfies \eqref{zero property fi}, \eqref{first property fi},\eqref{second property fi}, \eqref{fourth property fi} and \eqref{fifth property fi}. If, in addition, $E\not\subseteq L_1,$ then $\varphi_{cut}$ also satisfies \eqref{third property fi}.
\end{lem}

\begin{proof}
We prove that the function $s\to\frac1s||\sigma_sx^*||_E$ is decreasing. Let $s_2>s_1.$ We have $s_2=s_3s_1$ and $s_3>1.$ Therefore,
$$\frac1{s_2}||\sigma_{s_3}(\sigma_{s_1}(x^*))||_E\leq\frac{||\sigma_{s_3}||_{E\to E}}{s_2}||\sigma_{s_1}(x^*)||_E\leq\frac1{s_1}||\sigma_{s_1}(x^*)||_E,$$
since $||\sigma_{s_3}||_{E\to E}\leq s_3.$ It follows immediately that the limit in \eqref{definition of fi} exists.

\eqref{zero property fi} Trivial.

\eqref{first property fi} This follows from the fact that $||\sigma_s(x^*)||_E\leq s||x||_E.$

\eqref{second property fi} Since $y\prec\prec x,$ it follows that $\sigma_s(y^*)\prec\prec\sigma_s(x^*).$ Since $E$ is strictly symmetric, it follows that $||\sigma_s(y^*)||_E\leq||\sigma_s(x^*)||_E.$ Therefore,
$$\varphi(y)=\lim_{s\to\infty}\frac1s||\sigma_s(y^*)||_E\leq\lim_{s\to\infty}\frac1s||\sigma_s(x^*)||_E=\varphi(x).$$

\eqref{third property fi} Applying the semigroup property of the dilation operators $\sigma_{\tau},$
$$\lim_{s\rightarrow\infty}\frac1s||\sigma_s(\sigma_{\tau}(x^*))||_E=\tau\lim_{\tau\rightarrow\infty}\frac1{s\tau}||\sigma_{s\tau}(x^*)||_E=\tau\varphi(x).$$

\eqref{fourth property fi} By triangle inequality,
$$|\,||\sigma_s(x^*)||_E-||\sigma_s(y^*)||_E\,|\leq||\sigma_s(x^*-y^*)||_E.$$
Using \eqref{second property of rearrangement} and Lemma \ref{dilation property} one can obtain
$\sigma_s(x^*-y^*)\prec\prec\sigma_s((x-y)^*).$ Since $E$ is strictly symmetric,
$$|\,||\sigma_s(x^*)||_E-||\sigma_s(y^*)||_E\,|\leq||\sigma_s((x-y)^*)||_E.$$
Now, one can divide by $s$ and let $s\to\infty.$ Therefore,
$$|\varphi(x)-\varphi(y)|\leq\varphi(x-y)\leq||x-y||_E.$$
\eqref{fifth property fi} It follows from \eqref{first property of rearrangement} and Lemma \ref{dilation property}
that $\sigma_s((x+y)^*)\prec\prec\sigma_s(x^*)+\sigma_s(y^*).$ Therefore,
$$\varphi(x+y)=\lim_{s\to\infty}\frac1s||\sigma_s((x+y)^*)||_E
\leq\lim_{s\to\infty}\frac1s(||\sigma_s(x^*)||_E+||\sigma_s(y^*)||_E)=\varphi(x)+\varphi(y).$$

Existence and properties \eqref{zero property fi}-\eqref{fifth
property fi} of $\varphi_{fin}$ can be proved in a similar way.
Existence and properties \eqref{zero property fi}, \eqref{first
property fi},\eqref{second property fi},\eqref{third property
fi},\eqref{fifth property fi} of $\varphi_{cut}$ can be proved in
a similar way. Let us prove \eqref{third property fi} for
$\varphi_{cut}.$

\eqref{third property fi} Assume $E\not\subset L_1.$ By Lemma \ref{fi bounded null} below,
$\varphi(x^*\chi_{[\tau^{-1},1]})=\varphi_{cut}(x^*\chi_{[\tau^{-1},1]})=0.$
Hence,
$$\varphi(x^*\chi_{[0,\tau^{-1}]})\leq\varphi(x^*\chi_{[0,1]})\leq\varphi(x^*\chi_{[0,\tau^{-1}]})+\varphi(x^*\chi_{[\tau^{-1},1]})=\varphi(x^*\chi_{[0,\tau^{-1}]}).$$
Therefore,
$$\varphi_{cut}(\sigma_{\tau}(x^*))=\varphi(\sigma_{\tau}(x^*\chi_{[0,\tau^{-1}]}))=\tau\varphi(x^*\chi_{[0,\tau^{-1}]})=\tau\varphi_{cut}(x).$$
\end{proof}

\begin{lem}\label{fi bounded null} If $E=E(0,1)$ be a symmetric space on $(0,1)$ and $x\in L_{\infty}\cap E,$ then $\varphi(x)=0.$ If $E=E(0,\infty)$ be a symmetric space on $(0,\infty)$ and $x\in L_{\infty}\cap E,$ then $\varphi_{fin}(x)=0.$ If $E=E(0,\infty)\not\subseteq L_1$ and $x\in E\cap L_{\infty},$ then $\varphi_{cut}(x)=0.$
\end{lem}
\begin{proof} Clearly, $\varphi(x)=\varphi(x^*\chi_{[0,1]})\leq||x||_{\infty}\varphi(\chi_{[0,1]})$ in the first case. Similarly, $\varphi_{fin}(x)\leq||x||_{\infty}\varphi_{fin}(\chi_{[0,1]})$ ($\varphi_{cut}(x)\leq||x||_{\infty}\varphi_{cut}(\chi_{[0,1]})$) in the second (third) case. It is clear that $\varphi(\chi_{[0,1]})=0$ ($\varphi_{fin}(\chi_{[0,1]})=0$) in the first (second) case. Also, $E\not\subset L_1$ implies that $||\chi_{[0,n]}||_E=o(n)$ and, therefore, $\varphi_{cut}(\chi_{[0,1]})=0.$ The assertion follows immediately.
\end{proof}

\begin{rem} If $E$ is a separable symmetric space, then $E\cap L_{\infty}$ is a dense subset in $E$ (see \cite{KPS}). It follows now from the Lemmas \ref{fi bounded null} and \ref{general properties of fi} that functional $\varphi$ vanishes on every separable space $E.$
\end{rem}

\begin{lem}\label{main property of fi} Let $E$ be a strictly symmetric space. For functions $0\leq x_1,\ldots,x_k\in E$ and numbers $\lambda_1,\ldots,\lambda_k\geq0$
$$\varphi(\sum_{i=1}^k\lambda_ix_i)=\varphi(\sum_{i=1}^k\lambda_ix^*_i).$$
If $E=E(0,\infty),$ then the same is valid for $\varphi_{fin}.$ If, in addition, $E\not\subseteq L_1,$ then the same is valid for $\varphi_{cut}.$
\end{lem}

\begin{proof}
Applying the inequality \eqref{third property of rearrangement} $n$ times, we have for positive functions $x_1,\ldots,x_{2^n}$
$$(x_1^*+\ldots+x_{2^n}^*)\prec\prec2^n\sigma_{2^{-n}}(x_1+\ldots+x_{2^n}).$$
Therefore, by Lemma \ref{general properties of fi}\eqref{second property fi},
$$\varphi(x_1^*+\ldots+x_{2^n}^*)\leq\varphi(2^n\sigma_{2^{-n}}(x_1+\ldots+x_{2^n})).$$
By Lemma \ref{general properties of fi}\eqref{third property fi}, $2^k\varphi(\sigma_{2^{-k}}(z^*))=\varphi(z^*).$ Therefore,
$$\varphi(x_1^*+\ldots+x_{2^n}^*)\leq\varphi(x_1+\ldots+x_{2^n}).$$
Converse inequality follows trivially from \eqref{first property of rearrangement} and Lemma \ref{general properties of fi}\eqref{second property fi}.

The assertion of Lemma follows now from Lemma \ref{general properties of fi}\eqref{fourth property fi}.
\end{proof}

Note, that $y$ and $z$ in the Proposition below are arbitrary, that is $y,z$ do not necessary belong to $Q_+(x).$

\begin{prop}\label{nonlinearity of fi} Let $E$ be a symmetric space equipped with a Fatou norm. If $x\geq0\in E,$ then in each of the following cases there exists a decomposition $x=y+z,$ such that $y,z\geq 0$ and such that the following assertions hold.
\begin{enumerate}[i)]
\item\label{l61} If $E=E(0,1),$ then $\varphi(x)=\varphi(y)=\varphi(z).$
\item\label{l62} If $E=E(0,\infty)$ and $\varphi_{cut}(x)=0,$ then $\varphi(x)=\varphi(y)=\varphi(z).$
\item\label{l63} If $E=E(0,\infty),$ then $\varphi_{fin}(x)=\varphi_{fin}(y)=\varphi_{fin}(z).$
\item\label{l64} If $E=E(0,\infty),$ then $\varphi_{cut}(x)=\varphi_{cut}(y)=\varphi_{cut}(z).$
\end{enumerate}
\end{prop}

\begin{proof} We will prove only the first assertion. The proofs of the third and fourth assertions are exactly the same. The proof of the second assertion requires replacement of the interval $[\frac1m,\frac1n]$ with the interval $[n,m].$

We may assume that $x=x^*.$ Fix $n\in N.$ The sequence $\sigma_n(x\chi_{[\frac1m,\frac1n]})$ converges to $\sigma_n(x\chi_{[0,\frac1n]})$ almost everywhere when $m\to\infty.$

By the definition of Fatou norm,
$$||\sigma_n(x\chi_{[\frac1m,\frac1n]})||_E\to_{m}||\sigma_n(x\chi_{[0,\frac1n]})||_E.$$
For each $n\in N,$ one can select $f(n)>n,$ such that
$$||\sigma_n(x\chi_{[\frac1{f(n)},\frac1n]})||_E\geq(1-\frac1n)||\sigma_n(x\chi_{[0,\frac1n]})||_E.$$
Fix some $n_0$ and set $n_k=f^k(n_0),$ $k\in N.$ Here, $f^k=f\circ\ldots\circ f$ ($k$ times). Define
$$y=\sum_{k=0}^{\infty}x\chi_{[\frac1{n_{2k+1}},\frac1{n_{2k}}]},$$
$$z=\sum_{k=1}^{\infty}x\chi_{[\frac1{n_{2k}},\frac1{n_{2k-1}}]}.$$
It is clear, that
\begin{equation}\label{y estimate}
\frac1{n_{2k}}||\sigma_{n_{2k}}(y^*)||_E\geq\frac1{n_{2k}}||\sigma_{n_{2k}}(y)||_E\geq\frac1{n_{2k}}||\sigma_{n_{2k}}(x\chi_{[\frac1{n_{2k+1}},\frac1{n_{2k}}]})||_E.
\end{equation}
By definition of $n_k,$
\begin{equation}\label{from nk def}
\frac1{n_{2k}}||\sigma_{n_{2k}}(x\chi_{[\frac1{n_{2k+1}},\frac1{n_{2k}}]})||_E\geq\frac1{n_{2k}}(1-\frac1{n_{2k}})||\sigma_{n_{2k}}(x\chi_{[0,\frac1{n_{2k}}]})||_E.
\end{equation}
It follows from \eqref{y estimate} and \eqref{from nk def} that
\begin{equation}\label{vspom nerav}
\frac1{n_{2k}}||\sigma_{n_{2k}}(y^*)||_E\geq(1-\frac1{n_{2k}})\frac1{n_{2k}}||\sigma_{n_{2k}}(x\chi_{[0,\frac1{n_{2k}}]})||_E\geq(1-\frac1{n_{2k}})\varphi(x\chi_{[0,\frac1{n_{2k}}]}).
\end{equation}
By Lemma \ref{fi bounded null}, $\varphi(x\chi_{[\frac1{n_{2k}},1]})=0.$ Since $\varphi$ is convex, then
\begin{equation}\label{vspom rav}
\varphi(x\chi_{[0,\frac1{n_{2k}}]})\leq\varphi(x)\leq\varphi(x\chi_{[0,\frac1{n_{2k}}]})+\varphi(x\chi_{[\frac1{n_{2k}},1]})=\varphi(x\chi_{[0,\frac1{n_{2k}}]}).
\end{equation}
It follows from \eqref{vspom nerav} and \eqref{vspom rav} that
$$\frac1{n_{2k}}||\sigma_{n_{2k}}(y^*)||_E\geq(1-\frac1{n_{2k}})\varphi(x).$$
Passing to the limit, we obtain $\varphi(y)\geq\varphi(x).$ The converse inequality is obvious.

Hence, $\varphi(y)=\varphi(x)=\varphi(z),$ and this completes proof of the Proposition.
\end{proof}

\begin{lem}\label{fi constancy on qprime} If space $E$ is strictly symmetric, then $\varphi(y)=\varphi(x)$ for every  $y\in Q'(x).$ If, in addition, $E=E(0,\infty),$ then $\varphi_{fin}(y)=\varphi_{fin}(x)$ for every $y\in Q'(x).$ If $E\not\subseteq L_1,$ then $\varphi_{cut}(y)=\varphi_{cut}(x)$ for every $y\in Q'(x).$
\end{lem}
\begin{proof} Let
$$z=\sum_{i=1}^s\lambda_ix_i,$$
where $\lambda_i\geq0,$ $\sum_{i=1}^s\lambda_i=1,$ $x_i\geq0$ and $x_i^*=x.$ By Lemma \ref{main property of fi}, we obtain $\varphi(z)=\varphi(x).$ However, $y\in Q'(x)$ can be approximated by such $z.$ Since $\varphi$ is continuous in strictly symmetric spaces, the lemma follows readily.

The proofs are the same in cases of $\varphi_{fin}$ and $\varphi_{cut}.$
\end{proof}

If $A$ is a convex set, then function $\theta:A\to R$ is called midpoint additive if and only if
$$\theta(\frac12(y_1+y_2))=\frac12(\theta(y_1)+\theta(y_2)),\quad y_1,y_2\in A.$$

\begin{prop}\label{linearity of fi} Let $E$ be a strictly symmetric space and $x\in E.$ Then the following assertions hold.\begin{enumerate}[i)]
\item If $E=E(0,1),$ then $\varphi$ is midpoint additive on $Q_+(x).$
\item If $E=E(0,\infty),$ then $\varphi_{fin}$ is midpoint additive on $Q_+(x).$
\item If $E\not\subseteq L_1,$ then $\varphi_{cut}$ is midpoint additive on $Q_+(x).$
\end{enumerate}
\end{prop}

\begin{proof} We will only prove the first assertion. The proofs of the other two assertions are exactly the same.

Let $y\in \conv(\extr(\Omega_{+}(x))),$ so that
$$y=\sum_{i=1}^m\lambda_ix_i,$$
where $\lambda_i\geq0,$ $\sum_{i=1}^m\lambda_i=1,$ $x_i\geq0$ and $x_i^*=x^*\chi_{[0,\beta_i]}.$ Denote $z=\sum_{i=1}^m\lambda_ix^*\chi_{[0,\beta_i]}$ and $u=\sum_{\beta_i>0}\lambda_i x^*\chi_{[0,1]}.$ By Lemma \ref{main property of fi}, $\varphi(y)=\varphi(z).$

Since $|z-u|\in L_{\infty},$ then $\varphi(|u-z|)=0$ by Lemma \ref{fi bounded null}. By the triangle inequality,
$$\varphi(u)\leq\varphi(z)+\varphi(|u-z|)=\varphi(z)\leq\varphi(u)+\varphi(|u-z|)=\varphi(u).$$
Hence, $\varphi(y)=\varphi(u)=(\sum_{\beta_i>0}\lambda_i)\varphi(x).$ It is clear that last expression is midpoint additive on the set $\conv(\extr(\Omega_{+}(x))).$ By Lemma \ref{general properties of fi}, the functional $\varphi$ is continuous on $Q_{+}(x).$ Hence, it is midpoint additive on the set $Q_+(x).$
\end{proof}

\begin{prop}\label{linearity of fi infinite} Let $E=E(0,\infty)$ be a symmetric space on semi-axis equipped with a Fatou norm. Suppose that $E\not\subseteq L_1$ and $x\in E.$ If $\Omega_+(x)=Q_+(x),$ then $\varphi$ is midpoint additive on $\Omega_+(x).$
\end{prop}

\begin{proof} It follows from the Proposition \ref{linearity of fi} that $\varphi_{cut}$ is midpoint additive on $Q_+(x).$ By assumption, $\Omega_+(x)=Q_+(x).$ Hence, $\varphi_{cut}$ is midpoint additive on $\Omega_+(x).$ It follows now from Proposition \ref{nonlinearity of fi} \eqref{l64} that $\varphi_{cut}(x)=0.$ This assertion and Lemma 2 imply that $\varphi(x^*\chi_{[0,\beta]})=0$ for every finite $\beta.$

Let $y\in\conv(\extr(\Omega_+(x))).$ Hence,
$$y=\sum_{i=1}^m\lambda_ix_i,$$
where $\lambda_i\geq0,$ $\sum_{i=1}^m\lambda_i=1,$ $x_i\geq0$ and $x_i^*=x^*\chi_{[0,\beta_i]}.$ By convexity of $\varphi,$
$$\varphi(y)\leq\varphi(\sum_{\beta_i\in[0,\infty)}\lambda_ix_i)+\varphi(\sum_{\beta_i=\infty}\lambda_ix_i).$$
However,
$$0\leq\varphi(\sum_{\beta_i\in[0,\infty)}\lambda_ix_i)\leq\sum_{\beta_i\in[0,\infty)}\lambda_i\varphi(x^*\chi_{[0,\beta_i]})=0.$$
It then follows that
$$\varphi(y)\leq\varphi(\sum_{\beta_i=\infty}\lambda_ix_i).$$
The converse inequality is obvious. By Lemma \ref{main property of fi},
$$\varphi(y)=\varphi(\sum_{\beta_i=\infty}\lambda_ix_i)=\varphi(\sum_{\beta_i=\infty}\lambda_ix_i^*)=(\sum_{\beta_i=\infty}\lambda_i)\varphi(x).$$
Clearly, the last expression is midpoint additive on $\conv(\extr(\Omega_+(x))).$ Hence, the functional $\varphi$ is midpoint additive on $Q_+(x)=\Omega_+(x).$
\end{proof}

\begin{lem}\label{null fi for finite support} Let $E=E(0,\infty)$ be a strictly symmetric space on $(0,\infty)$ and $x\in E.$ Suppose, that $E\not\subseteq L_1.$ If $P(x|\mathcal{A})\in Q'(x)$ for every $\mathcal{A},$ then $\varphi_{cut}(x)=0.$
\end{lem}

\begin{proof} Suppose that $x=x^*.$ Set $\mathcal{A}=\{[0,1]\}$ and $y=P(x|\mathcal{A})\in E\cap L_{\infty}.$ By assumption, $y\in Q'(x).$ By Lemma \ref{fi constancy on qprime} and Lemma \ref{fi bounded null}, $\varphi_{cut}(x)=\varphi_{cut}(y)=0.$
\end{proof}

\begin{lem}\label{if linfty in e} Let $E$ and $x$ be as in Lemma \ref{null fi for finite support}. If $L_{\infty}\subseteq E,$ then $\varphi(x)=0.$
\end{lem}
\begin{proof} Due to the choice of $E,$ we have $1\in E.$ However, $\sigma_{\tau}(1)=1$ implies $\varphi(1)=0.$ Thus, for every $z\in E\cap L_{\infty},$ we have $\varphi(z)=0.$ However, for every $x\in E,$ we have $\varphi(x^*\chi_{[0,1]})=0$ due to Lemma \ref{null fi for finite support}. Hence,
$$0\leq\varphi(x)=\varphi(x^*)\leq\varphi(x^*\chi_{[0,1]})+\varphi(x^*\chi_{[1,\infty)})=0+0=0.$$
\end{proof}

\begin{lem}\label{fi exact} Let $E$ and $x$ be as in Lemma \ref{null fi for finite support}. If $y\in E\cap L_{\infty}$ and if
$$\omega(y)=\limsup_{t\to\infty}\frac{\int_0^ty^*(s)ds}{\int_0^tx^*(s)ds},$$
then $\varphi(y)=\omega(y)\varphi(x).$
\end{lem}

\begin{proof} Fix $\varepsilon>0.$ There exists $T>0,$ such that for every $t>T,$
$$\int_0^ty^*(s)\leq(\omega(y)+\varepsilon)\int_0^tx^*(s)ds.$$
It then follows that $y\prec\prec(\omega(y)+\varepsilon)(x^*+C\chi_{[0,T]})$ for some constant $C.$ By Lemma \ref{general properties of fi} \eqref{second property fi}, $\varphi(y)\leq(\omega(y)+\varepsilon)\varphi(x^*+C\chi_{[0,T]}).$ By Lemma \ref{fi bounded null}, $\varphi(C\chi_{[0,T]})=0$ and, therefore, $\varphi(x^*+C\chi_{[0,T]})=\varphi(x).$ Hence $\varphi(y)\leq\omega(y)\varphi(x).$

Now, fix $\omega<\omega(y).$ There exists a sequence $t_k\to\infty,$ such that
$$\int_0^{t_k}y^*(s)ds\geq\omega\int_0^{t_k}x^*(s)ds.$$
Without loss of generality, $t_0=0.$ Set $u=P(x^*|\mathcal{A}),$ where $\mathcal{A}=\{[t_k,t_{k+1})\}.$ It then follows that $\omega u\prec\prec y$ and $\omega\varphi(u)\leq\varphi(y).$ However, $u\in Q'(x)$ and $\varphi(u)=\varphi(x)$ due to Lemma \ref{fi constancy on qprime}. Hence $\omega(y)\varphi(x)\leq\varphi(y).$
\end{proof}

\begin{prop}\label{cut function lemma} Let $E=E(0,\infty)$ be a symmetric space on the semi-axis and let $x\in E.$ If $\varphi(x)=0,$ then, $x\chi_A\in Q'(x)$ for every Lebesgue measurable subset $A\subseteq(0,\infty).$
\end{prop}

\begin{proof} Let $[0,\infty)=B\cup C,$ where $B,C$ are disjoint sets such that $m(B)=m(A)$ and $m(C)=\infty.$ Fix a partition $C=\cup_{i=1}^{n+1}C_i,$ where $m(C_i)=m(R_+\backslash A),$ $1\leq i\leq n.$ Let $\gamma:B\to A$ and $\gamma_i:C_i\to R_+\backslash A,$ $1\leq i\leq n$ be measure-preserving transformations. Define functions $x_n^i,$ $1\leq i\leq n$ by the following construction. Set $x_n^i\chi_B=x\circ\gamma,$ $x_n^i|_{C_i}=x\circ\gamma_i$ and $x_n^i|_{C_j}=0$ if $i\neq j.$ Clearly, $x_n^i\sim x$ and
$$||(x\chi_A)\circ\gamma-\frac1n\sum_{i=1}^nx_n^i||_E=\frac1n||\sigma_n(x\chi_{[0,\infty)\backslash A})||_E\leq\frac1n||\sigma_n(x^*)||_E\to 0.$$
Hence, $(x\chi_A)\circ\gamma\in Q'(x).$ Thus, $x\chi_A\in Q'(x).$
\end{proof}

\begin{cor}\label{main suff coroll} Let $E=E(0,\infty)$ be a symmetric space on semi-axis. If $\varphi(x)=0,$ then $y\chi_A\in Q'(x)$ for every $y\in Q'(x).$
\end{cor}
\begin{proof} It follows from assumption and Lemma \ref{fi constancy on qprime} that $\varphi(y)=\varphi(x)=0.$ Lemma \ref{cut function lemma} implies that $y\chi_{A}\in Q'(y).$  Since $y\chi_A\in Q'(y)$ and $y\in Q'(x),$ then Lemma \ref{transitivity lemma} implies $y\chi_A\in Q'(x).$
\end{proof}

An assertion somewhat similar to the lemma below is contained in \cite[Lemma 1.3]{BM}.

\begin{lem}\label{normality of qprime} Assume that $x\in E$ satisfies conditions of Proposition \ref{cut function lemma}. If $y\in Q'(x)$ and $0\leq z\leq y,$ then, $z\in Q'(x).$
\end{lem}

\begin{proof} Define sets $e_n^i,$ $i=1,\ldots,n$ by
$$e_n^i=\{t:\ \frac{i-1}ny(t)\leq z(t)\leq\frac{i}ny(t)\}.$$
Define functions $y_n^k,$ $k=1,\ldots,n$ as $y_n^k=y\sum_{k<(i+n)/2}\chi_{e_n^i}.$ By Corollary \ref{main suff coroll}, $y_n^k\in Q'(x).$ Put
$$s_n=\frac1n\sum_{k=1}^ny_n^k\in Q'(x).$$
Clearly,
$$|s_n(t)-(y(t)+z(t))/2|\leq\frac{2y(t)}n,\ \ \forall t\in e_n^i.$$
Hence, $s_n\to(y+z)/2$ by norm. Therefore, $(y+z)/2\in Q'(x).$ We can repeat this procedure $n$ times and obtain $2^{-n}((2^n-1)z+y)\in Q'(x).$ Therefore, $z\in Q'(x).$
\end{proof}

The following assertion seems to be known. We include details of proof for lack of a convenient reference.

\begin{lem}\label{transitivity lemma} Let $E$ be a symmetric space either on $(0,1)$ or $(0,\infty)$ and $x\in E.$ If $y\in Q'(z)$ and $z\in Q'(x),$ then $y\in Q'(x).$
\end{lem}

\begin{proof} Without loss of generality, $y=y^*,$ $z=z^*$ and $x=x^*.$ Let $y\in Q'(z).$ Hence, for every $\varepsilon>0,$ one can find $n\in N,$ $\lambda_i\in R_+$ and measurable functions $z_i\sim z,$ $i=1,\ldots,n,$ such that $\sum_{i=1}^n\lambda_i=1$ and
$$||y-\sum_{i=1}^n\lambda_iz_i||_E\leq\varepsilon.$$
One can find measure-preserving transformations $\gamma_i,$ such that
$$||z_i-z\circ\gamma_i||_{L_1\cap L_{\infty}}\leq\varepsilon.$$
Hence,
$$||y-\sum_{i=1}^n\lambda_iz\circ\gamma_i||_E\leq2\varepsilon.$$
However, $z\in Q'(x).$ Consequently, arguing in a similar way, one can find $m\in N,$ $\mu_j\in R_+$ and measure preserving transforamtions $\delta_j,$ $1\leq j\leq n$ such that $\sum_{j=1}^m\mu_j=1$ and
$$||z-\sum_{j=1}^m\mu_jx\circ\delta_j||_E\leq2\varepsilon.$$
Therefore,
$$||y-\sum_{i=1}^n\sum_{j=1}^m\lambda_i\mu_jx\circ\gamma_i\circ\delta_j||\leq4\varepsilon$$
and this suffices to complete the proof.
\end{proof}

\begin{rem} The collection of sets $\{Q(x),\ x\in E\}$ also satisfies the transitivity property expressed in Lemma \ref{transitivity lemma}. We do not know whether this is the case for the collection $\{Q_+(x),\ x\in E\}.$
\end{rem}

\section{Main results}

The implication $(ii)\Rightarrow(i)$ in the following theorem is almost verbatim repetition of the argument given in \cite[Lemma 3.1]{BM} for the case of finite measure. For convenience of the reader, we present here a proof of the most important case.

\begin{thm}\label{expectation lemma} $(a)$ Let $E$ be a fully symmetric space and $x\in E.$ If $E=E(0,1)$ or $E=E(0,\infty)$ and $E\not\subseteq L_1,$ then the following conditions are equivalent.
\begin{enumerate}[i)]
\item $P(x|\mathcal{A})\in Q'(x)$ for every $\mathcal{A}\in\mathfrak{A}.$
\item $\varphi(x)=0.$
\end{enumerate}
$(b)$ If $E=E(0,\infty)$ and $E\subseteq L_1,$ then the following conditions are equivalent.
\begin{enumerate}[i)]
\item $P(x|\mathcal{A})\in Q'(x)$ for every $\mathcal{A}\in\mathfrak{A}.$
\item $\varphi_{fin}(x)=0.$
\end{enumerate}

\end{thm}

\begin{proof} $(a)$ $(i)\Rightarrow(ii)$ Let $E=E(0,1)$ and $x=x^*.$ Set $\mathcal{A}=\{[0,1]\}$ and $y=P(x|\mathcal{A}).$ By assumption, $y\in Q'(x).$ By Lemma \ref{fi constancy on qprime} and Lemma \ref{fi bounded null}, $\varphi(x)=\varphi(y)=0.$

Let $E=E(0,\infty)$ and $L_{\infty}\subseteq E\not\subseteq L_1.$ The assertion is proved in Lemma \ref{if linfty in e}.

Let $E=E(0,\infty)$ and $L_{\infty}\not\subseteq E\not\subseteq L_1.$ Suppose that $x=x^*$ and $\varphi(x)>0.$ Set $\mathcal{B}=\{[0,1]\},$ $\psi'=P(x|\mathcal{B})$ and $\psi(t)=\int_0^t\psi'(s)ds.$ By Lemma \ref{fi constancy on qprime}, $\varphi(\psi')=\varphi(x).$

Let $y\in E\cap L_{\infty}.$ It follows from Lemma \ref{fi exact}, that $\omega(y)\varphi(x)<\infty.$ Since $\varphi(x)>0,$ then $\omega(y)<\infty.$ Therefore, $y\in M_{\psi}.$ Hence, $E\cap L_{\infty}\subseteq M_{\psi}.$ Since $E$ is fully symmetric and $\psi'\in E\cap L_{\infty},$ then $M_{\psi}\in E\cap L_{\infty}.$ Therefore, $E\cap L_{\infty}=M_{\psi}.$

If $u=2\sigma_{\frac12}\psi',$ then $\varphi(u)=\varphi(\psi')$ by Lemma \ref{general properties of fi}\eqref{fourth property fi}. Hence $\omega(u)\varphi(x)=\varphi(x)$ and $\omega(u)=1.$ However,
$$\omega(u)=\limsup_{t\to\infty}\frac{\int_0^t2x(2s)ds}{\int_0^tx(s)ds}=\limsup_{t\to\infty}\frac{\psi(2t)}{\psi(t)}.$$
Thus,
\begin{equation}\label{russu condition}
\lim_{t\to\infty}\frac{\psi(2t)}{\psi(t)}=1.
\end{equation}
Let $G$ be the set defined by
$$G=\{y\in E:\ \exists C\ \sup_{t\geq1}\frac{y^*(t)}{\psi'(Ct)}<\infty\}.$$
Note, that our set $G$ differs from the one introduced in \cite{KPS}. If $y_1,y_2\in G,$ then $y^*_i(t)\leq C_i\psi(Ct)$ for $t\geq\frac12.$ It then follows
$$(y_1+y_2)^*(t)\leq y^*_1(\frac{t}2)+y^*_2(\frac{t}2)\leq (C_1+C_2)\psi'(\frac{C}2t).$$
In particular, $G$ is a linear set and $\conv(\{y^*=x^*\})\subseteq G.$ If the condition \eqref{russu condition} holds, then there exists a sequence $t_k,$ such that $t_0=0,$ $t_1=1$ and for every $k$
$$\frac{\psi(t_{k+1})-\psi(t_k)}{t_{k+1}-t_k}\geq\frac23\frac{\psi(\frac12t_{k+1})}{t_{k+1}}.$$
Set $\mathcal{A}=\{[t_k,t_{k+1}]\}$ and $z=P(x|\mathcal{A}).$ It follows from the construction given in \cite{KPS} that  $||(z-y)\chi_{[\frac12t_{k},t_k]}||_{M_{\psi}}\geq\frac14$ for every $y\in G$ and every sufficiently large $k.$ However, $||(y-z)\chi_{[\frac12t_{k},t_k]}||_{L_{\infty}}\to 0.$ Since $M_{\psi}=E\cap L_{\infty},$ then $||(z-y)\chi_{[\frac12t_{k},t_k]}||_{E}\geq\frac14$ for sufficiently large $k.$ In particular, $||y-z||_E\geq\frac14.$ Hence, $\dist_E(z,G)\geq\frac14$ and $\dist_E(z,Q'(x))\geq\frac14.$ This contradicts the assumption that $P(x|\mathcal{A})\in Q'(x).$

$(a)$ $(ii)\Rightarrow(i)$ Let $E=E(0,1)$ or $E=E(0,\infty)\not\subseteq L_1.$ We will prove the assertion for the case, when $\mathcal{A}=\{[0,1]\}.$ The general proof is similar. Without loss of generality, $x$ decreases on $[0,1].$ Define functions $x_n^i,$ $i=0,\ldots,n-1$ such that $(i)$ $x_n^i=x$ outside $(0,1)$ and $(ii)$ $x_n^i(t)=x(t+\frac{i}n(\mmod 1))$ if $t\in(0,1).$ Set $x_n(t)=x(t-\frac{i}n)$ if $\frac{i}n\leq t\leq\frac{i+1}n,$ $0\leq i\leq n-1$ and $x_n(t)=0$ if $t\geq 1.$ Clearly, $x_n^i\sim x$ and $(x_n)^*\leq\sigma_n(x^*).$

We will show that
$$\int_0^1x(s)ds-\frac1n\sum_{i=0}^{n-1}x(t+\frac{i}n(\mmod 1))\leq\int_0^{\frac1n}x(s)ds$$
and
$$\int_0^1x(s)ds-\frac1n\sum_{i=0}^{n-1}x(t+\frac{i}n(\mmod1))\geq-\frac1nx_n(t).$$
We will prove only the first inequality. The proof of the second one is identical.

Without loss of generality, $t\in[0,\frac1n].$ Clearly,
$$\frac1nx(t+\frac{i}n)\geq\int_{\frac{i+1}n}^{\frac{i+2}n}x(s)ds$$
for $i=0,\ldots,n-2.$ Hence,
$$\int_0^1x(s)ds-\frac1n\sum_{i=0}^{n-1}x(t+\frac{i}n)=\int_0^{\frac1n}x(s)ds-\frac1n x(t+\frac{n-1}n)-$$
$$-\sum_{i=0}^{n-2}(\frac1nx(t+\frac{i}n)-\int_{\frac{i+1}n}^{\frac{i+2}n}x(s)ds)\leq\int_0^{\frac1n}x(s)ds.$$

Therefore,
$$|\int_0^1x(s)ds-\frac1n\sum_{i=0}^{n-1}x_n^i(t)|\leq\frac1nz_n(t),\quad t\in[0,1],$$
where $z_n=x_n+(\int_0^1x_n(s)ds)\chi_{[0,1]}.$ Obviously, $z_n\prec\prec2x_n\leq2\sigma_n(x^*)$ and, therefore, $||z_n||_E\leq2||\sigma_n(x^*)||_E.$

It then follows that
$$||P(x|\mathcal{A})-\frac1n\sum_{i=0}^{n-1}x_n^i||_E\leq\frac2n||\sigma_n(x^*)||_E\to0.$$

$(b)$ $(i)\Rightarrow(ii)$ Let $E=E(0,\infty)$ and $E\subset L_1.$ Set $A=\{s:x(s)\geq 1\}$ and $\mathcal{A}=\{A\}.$ Set $y=P(x|\mathcal{A})\in E\cap L_{\infty}.$ Lemma \ref{fi bounded null} implies that $\varphi(y)=0.$ By assumption, $y\in Q'(x).$ By Lemma \ref{fi constancy on qprime}, $\varphi(x)=\varphi(y)=0.$

$(b)$ $(ii)\Rightarrow(i)$ The assertion follows from Theorem \ref{qprime integrable}.
\end{proof}

The following proposition is the core technical result of the article. In case of the interval $(0,1)$ it may be found in \cite[Lemma 3.2]{BM}. However, our proof is more general, simpler and shorter.

We consider functions of the form
\begin{equation}\label{elementary}
x=\sum_{i\in Z}x_i\chi_{[a_{i-1},a_i]},\ \ y=\sum_{i\in Z}y_i\chi_{[a_{i-1},a_i]},
\end{equation}
where $\{a_i\}_{i\in Z}$ is an increasing sequence (possibly finite or one-sidedly infinite).

\begin{prop}\label{Braverman-Mekler lemma} Let $y=y^*$ and $x=x^*$ be functions of the form \eqref{elementary} either on $(0,1)$ or on $(0,\infty).$ If $y\prec\prec x,$ then there exists a countable collection $\{\Delta_k\}_{k\in\mathcal{K}}$ of disjoint sets, where $\Delta_k=I_k\cup J_k$ with intervals $I_k$ and $J_k$ of finite measure, such that
\begin{enumerate}[i)]
\item\label{first item of bm} The functions $x$ and $y$ are constant on the intervals $I_k$ and $J_k$ and the interval $I_k$ lies to the left of $J_k,$ $k\in\mathcal{K}.$
\item\label{second item of bm} $y|_{\Delta_k}\prec x|_{\Delta_k},$ $k\in\mathcal{K}.$
\item\label{third item of bm} $y(t)\leq x(t)$ if $t\notin\cup_{k\in\mathcal{K}}\Delta_k.$
\end{enumerate}
If, in addition, $x$ and $y$ are functions on $(0,1)$ and $\int_0^1y(s)ds=\int_0^1x(s)ds,$ then $y(t)=x(t)$ if $t\notin\cup_{k\in\mathcal{K}}\Delta_k.$
\end{prop}
\begin{proof}
There exists a subsequence $\{a_{m_i}\}_{i\in\mathcal{I}}$ (possibly finite or one-sidedly infinite) such that  $\{x<y\}=\cup_{\in\mathcal{I}}[a_{m_i-1},a_{m_i}].$  Since $y\prec\prec x,$ we have
$$\int_0^t(x-y)_+(s)ds-\int_0^t(y-x)_+(s)ds=\int_0^tx(s)ds-\int_0^ty(s)ds\geq0.$$
For each $i\in\mathcal{I},$ denote by $b_i$ the minimal $t>0,$ such that
$$\int_0^t(x-y)_+(s)ds=\int_0^{a_{m_i}}(y-x)_+(s)ds.$$
Clearly, for every $i\in\mathcal{I},$
$$\int_0^{a_{m_i-1}}(x-y)_+(s)ds=\int_0^{a_{m_i}}(x-y)_+(s)ds\geq\int_0^{a_{m_i}}(y-x)_+(s)ds.$$
Hence, $b_i\leq a_{m_i-1}.$ For each $i\in\mathcal{I},$ the set $[b_{i-1},b_i]\cap\{x>y\}$ is a finite union $\cup_{j=1}^{n_i}I_i^j$ of disjoint intervals on which each of $x$ and $y$ is finite. By the definition of $b_i,$ we have
$$\int_{a_{m_i-1}}^{a_{m_i}}(y-x)_+(s)ds=\int_{b_{i-1}}^{b_i}(x-y)_+(s)ds=\sum_{j=1}^{n_i}\int_{I_i^j}(x-y)_+(s)ds.$$
Set $\mathcal{K}=\{(i,j):1\leq j\leq n_i, i\in\mathcal{I}\}.$ If $k=(i,j)\in\mathcal{K},$ set $I_k=I_i^j$ and
$$J_k=J_i^j=[a_{m_i-1}+(y_{m_i}-x_{m_i})^{-1}c_i^{j-1},a_{m_i-1}+(y_{m_i}-x_{m_i})^{-1}c_i^j],$$
where
$$c_i^j=\sum_{l=1}^j\int_{I_i^l}(x-y)_+(s)ds,\quad i\in\mathcal{I},0\leq j\leq n_i.$$
Using the fact that $x$ and $y$ are constant on the interval $[a_{m_i-1},a_{m_i}],$ we obtain $J_k\subset[a_{m_i-1},a_{m_i}]$ and $\cup_{j=1}^{n_i}J_i^j=[a_{m_i-1},a_{m_i}].$

\eqref{first item of bm} Both $x$ and $y$ are constant on $I_k$ and $J_k,$ $k\in\mathcal{K}.$ Since $b_{i}\leq a_{m_i-1}$ for each $i\in\mathcal{I},$ then $I_k$ lies to the left of $J_k$ for $k\in\mathcal{K}.$

It then follows from \eqref{first item of bm}, that
\begin{equation}\label{integrals coincide}
\int_{I_k}(x-y)_+(s)ds=\int_{J_k}(y-x)_+(s)ds,\quad k\in\mathcal{K}.
\end{equation}

\eqref{second item of bm} Since $x|_{I_k}\geq y|_{I_k}$ and $x|_{J_k}\leq y|_{J_k}$ for all $k\in\mathcal{K},$ then the assertion follows directly from \eqref{first item of bm} and \eqref{integrals coincide}.

\eqref{third item of bm} The set $\{y>x\}=\cup_{i\in\mathcal{I}}\cup_{j=1}^{n_i}J_i^j\subseteq\cup_{k\in\mathcal{K}}\Delta_k.$

The last assertion is immediate.
\end{proof}

\begin{cor}\label{tech corollary} Let $E$ be a fully symmetric space either on the interval $(0,1)$ or on the semi-axis. If $x,$ $y$ and $\mathcal{B}=\{\Delta_k\}_{k\in\mathcal{K}}$ are as in Proposition \ref{Braverman-Mekler lemma} and $y(t)=x(t)$ if $t\notin\cup_k\Delta_k,$ then $y$ can be arbitrary well approximated in the norm of $E$ by convex combinations of functions of the form $P(x|\mathcal{A}),$ $\mathcal{A}\in\mathfrak{A}.$
\end{cor}

\begin{proof} Set $\lambda_k=(y|_{I_k}-y|_{J_k})/(x|_{I_k}-x|_{J_k}),$ $k\in\mathcal{K}.$ Since $y|_{\Delta_k}\prec x|_{\Delta_k},$ it is not difficult to verify that $\lambda_k\in[0,1],$ $k\in\mathcal{K}.$ Further, a simple calculation shows that $y=(1-\lambda_k)P(x|\mathcal{B})+\lambda_kx$ on $\Delta_k,$ $k\in\mathcal{K}.$

As well-known, every $[0,1]-$valued sequence can be uniformly approximated by convex combinations of $\{0,1\}-$valued sequences.

Fix $\varepsilon>0.$ There exists $\mu\in l_{\infty}(\mathcal{K})$ with $\mu=\sum_{i=1}^n\theta_i\chi_{D_i}$ for some $n\in N,$ $0\leq\theta_i\in R$ and $D_i\subseteq\mathcal{K}$ such that $\sum_{i=1}^n\theta_i=1$ and $||\lambda-\mu||_{\infty}\leq\varepsilon.$
Set $z=(1-\mu_k)P(x|\mathcal{B})+\mu_k x$ on $\Delta_k,$ $k\in\mathcal{K}$ and $z=x$ outside $\cup_{k\in\mathcal{K}}\Delta_k.$ It is clear that $|y-z|\chi_{\Delta_k}=|\lambda_k-\mu_k|\,|x-P(x|\mathcal{B})|\chi_{\Delta_k},$ $k\in\mathcal{K}$ and $|y-z|=\sum_{k\in\mathcal{K}}|y-z|\chi_{\Delta_k}\leq2\varepsilon(x+P(x|\mathcal{B})).$ Therefore, $||y-z||_E\leq2\varepsilon||x||_E.$

Set $F_i=\cup_{k\in D_i}\Delta_k$ and  $\mathcal{A}_i=\{\Delta_k\}_{k\notin D_i}\in\mathfrak{A},$ $1\leq i\leq n.$ It is then clear that
$$z=\sum_{i=1}^n\theta_i((1-\chi_{F_i})P(x|\mathcal{B})+\chi_{F_i}x)=\sum_{i=1}^n\theta_iP(x|\mathcal{A}_i).$$
\end{proof}

\subsection{The case that $E\subseteq L_1$}

\begin{thm}\label{qprime finite} Let $E=E(0,1)$ be a fully symmetric space on the interval $(0,1).$
If $x\in E,$ then the following statements are equivalent.
\begin{enumerate}[i)]
\item $\Omega'(x)=Q'(x).$
\item $\varphi(x)=0.$
\end{enumerate}
\end{thm}

\begin{proof} $(i)\Rightarrow(ii)$ Suppose that $Q'(x)=\Omega'(x).$ Set $\mathcal{A}=\{[0,1]\}$ and
$y=P(x|\mathcal{A}).$ Clearly, $y\in\Omega'(x)=Q'(x).$ Lemma
\ref{fi constancy on qprime} implies that $\varphi(x)=\varphi(y).$
Lemma \ref{fi bounded null} implies $\varphi(y)=0.$ The assertion
is proved.

$(ii)\Rightarrow(i)$ Let $x=x^*$ and $0\leq y\in\Omega'(x).$ In
this case, $y=y^*\circ\gamma$ for some measure-preserving
transformation $\gamma$ (see \cite{Ryff3} or \cite[Theorem 7.5,
p.82]{BS}). Without loss of generality, we may assume that
$y=y^*.$ Fix $\varepsilon>0.$ Set
$$s_n(\varepsilon)=\inf\{s: y(s)\leq y(1)+n\varepsilon\},\quad n\in N.$$
Let $\mathcal{A}_{\varepsilon}$ be the partition, determined by the points  $s_n(\varepsilon),$ $n\in N.$ Set $u=P(y|\mathcal{A}_{\varepsilon})$ and $z=P(x|\mathcal{A}_{\varepsilon}).$ The functions $u$ and $z$ satisfy the condition $u\prec z$ and are of the form given in \eqref{elementary}.

By Lemma \ref{general properties of fi}\eqref{second property fi}, $\varphi(z)\leq\varphi(x)=0.$ By Theorem \ref{expectation lemma}, $P(z|\mathcal{A})\in Q'(z)$ for every $\mathcal{A}\in\mathfrak{A}$ It follows now from Corollary \ref{tech corollary} that $u\in Q'(z).$ However, $z\in Q'(x)$ by Theorem \ref{expectation lemma}. Therefore, by Lemma \ref{transitivity lemma}, $u\in Q'(x).$ However, $||y-u||_{L_{\infty}}\leq\varepsilon.$ Since $\varepsilon$ is arbitrary, $y\in Q'(x).$
\end{proof}

\begin{thm}\label{qplus finite} Let $E=E(0,1)$ be a fully symmetric space on the interval $(0,1).$ If $x\in E$ and $\varphi(x)=0,$ then $\Omega_+(x)=Q_+(x).$ If, in addition, the norm on $E$ is a Fatou norm, then converse assertion also holds.
\end{thm}
\begin{proof} Suppose that $\varphi(x)=0$ and let $y\in\Omega_+(x).$ Hence, there exists $s_0\in[0,1],$ such that $\int_0^{s_0}x^*(s)ds=\int_0^1y^*(s)ds.$ Set $z=x^*\chi_{[0,s_0]}.$ By Theorem \ref{qprime finite}, $y\in Q'(z).$ Hence, $y\in Q'(z)\subseteq Q_+(x).$

By Proposition \ref{nonlinearity of fi}, there exist $0\leq y,z\in E,$ such that $x=y+z$ and $\varphi(x)=\varphi(y)=\varphi(z).$ By Proposition \ref{linearity of fi}, $\varphi(x)=\varphi(y)+\varphi(z).$ Consequently, $\varphi(x)=0.$
\end{proof}

Now, consider the case that $E=E(0,\infty).$

\begin{thm}\label{qprime integrable} Let $E=E(0,\infty)$ be a fully symmetric space on semi-axis. If $E\subseteq L_1$ and $x\in E,$ then the following assertions are equivalent.
\begin{enumerate}[i)]
\item $\Omega'(x)=Q'(x).$
\item $\varphi_{fin}(x)=0.$
\end{enumerate}
\end{thm}

\begin{proof} $(i)\Rightarrow(ii)$  Let $x=x^*$ and suppose that $Q'(x)=\Omega'(x).$ Set $\mathcal{A}=\{[0,1]\}$ and $y=P(x|\mathcal{A}).$ Clearly, $y\in\Omega'(x)=Q'(x).$ Lemma \ref{fi constancy on qprime} implies that $\varphi(x)=\varphi(y).$ Lemma \ref{fi bounded null} implies $\varphi(y)=0.$ The assertion is proved.

$(ii)\Rightarrow(i)$ Let $x=x^*$ and $0\leq y\in\Omega'(x).$ It
follows from \cite[Lemma II.2.1]{KPS} that for every fixed
$\varepsilon>0$ there exists measure-preserving transformation
$\gamma$ such that $||y-y^*\circ\gamma||_E\leq\varepsilon.$
Without loss of generality, we may assume that $y=y^*.$ For every
$S>0,$
$$\frac1{\tau}||(\sigma_{\tau}x)\chi_{[0,S]}||_E\leq\frac{S}{\tau}||(\sigma_{\tau}x)\chi_{[0,1]}||_E\to 0.$$

$(a)$ Suppose first that $\supp(x)=\supp(y)=(0,\infty).$ Fix $\varepsilon>0.$ There exists $T,$ such that
$$||x\chi_{[T,\infty)}||_{L_1\cap L_{\infty}}\leq\varepsilon,\quad ||y\chi_{[T,\infty)}||_{L_1\cap L_{\infty}}\leq\varepsilon.$$
Clearly, $\int_0^Tx(s)ds<\int_0^{\infty}x(s)ds.$ Hence, there exists $S\geq T,$ such that $\int_0^Sy(s)ds=\int_0^Tx(s)ds.$ By Theorem \ref{qprime finite}, $y\chi_{[0,S]}\in Q'(x\chi_{[0,T]}).$ Hence, $y\in Q'(x)+y\chi_{(S,\infty)}-Q'(x\chi_{(T,\infty)})$ and, therefore, $\dist(y,Q'(x))\leq 2\varepsilon.$ Since $\varepsilon$ is arbitrary, $y\in Q'(x).$

$(b)$ Suppose now that $m(\supp(x))<\infty$ or $m(\supp(y))=0.$ Fix $z=z^*\in L_1\cap L_{\infty}$ with infinite support. It is clear that $(y+\varepsilon z)\in\Omega'(x+\varepsilon z),$ $\varepsilon>0.$ By assumption and Lemma \ref{fi bounded null}, $\varphi_{fin}(x+\varepsilon z)=0.$ Hence, using $(a)$ preceding, it follows that $(y+\varepsilon z)\in Q'(x+\varepsilon z)\subset Q'(x)+\varepsilon Q'(z).$ Hence, $\dist(y,Q'(x))\leq\varepsilon$ for every $\varepsilon>0$ and, therefore, $y\in Q'(x).$
\end{proof}

\begin{thm}\label{qplus integrable} Let $E=E(0,\infty)$ be a fully symmetric space on $(0,\infty)$ such that $E\subseteq L_1.$ If $0\leq x\in E$ and $\varphi_{fin}(x)=0,$ then $\Omega_+(x)=Q_+(x).$ If, in addition, the norm on $E$ is a Fatou norm, then converse assertion also holds.
\end{thm}

\begin{proof} Let $\varphi_{fin}(x)=0$ and $y\in\Omega_+(x).$ As in Theorem \ref{qprime integrable}, we may assume $y=y^*.$ Fix $\varepsilon>0.$ There exists $T>0$ such that $$||x\chi_{[T,\infty)}||_{L_1\cap L_{\infty}}\leq\varepsilon,\ ||y\chi_{[T,\infty)}||_{L_1\cap L_{\infty}}\leq\varepsilon.$$
Select $S\leq T$ such that
$$\int_0^Sx^*(s)ds=\int_0^Ty^*(s)ds.$$
Clearly, $y\chi_{[0,T]}\in\Omega'(x^*\chi_{[0,S]}).$ By Theorem \ref{qprime finite}, $y\chi_{[0,T]}\in Q'(x^*\chi_{[0,S]})\subseteq Q_+(x).$ Hence, $y\in Q_+(x).$

By Proposition \ref{nonlinearity of fi}, there exist $0\leq y,z\in E,$ such that $x=y+z$ and $\varphi_{fin}(x)=\varphi_{fin}(y)=\varphi_{fin}(z).$ By Proposition \ref{linearity of fi}, $\varphi_{fin}(x)=\varphi_{fin}(y)+\varphi_{fin}(z).$ Consequently, $\varphi_{fin}(x)=0.$
\end{proof}

\subsection{The case that $E\not\subseteq L_1$}

\begin{thm}\label{qprime nonintegrable} Let $E=E(0,\infty)$ be a fully symmetric space on the semi-axis and let $x\in E.$ If $\varphi(x)=0,$ then $\Omega_+(x)=Q'(x).$
\end{thm}

\begin{proof} Let us assume first that $y=y^*\in\Omega_+(x).$ Fix $\varepsilon>0.$ Set $t_n(\varepsilon)=1+n\varepsilon,$
$$s_n(\varepsilon)=\inf\{s: y(s)\leq y(1)+n\varepsilon\},$$
$$s_{-n}(\varepsilon)=\sup\{s: y(s)\geq y(1)-n\varepsilon\}.$$
Let $\mathcal{A}_{\varepsilon}$ be the partition, determined by the points  $s_{\pm n}(\varepsilon),$ $t_n(\varepsilon).$ Set $u=P(y|\mathcal{A}_{\varepsilon})$ and $z=P(x|\mathcal{A}_{\varepsilon}).$ The functions $u$ and $z$ satisfy the conditions $u\prec\prec z$ and \eqref{elementary}. Set
$$v=u\sum_{k\in\mathcal{K}}\chi_{\Delta_k}+z\chi_{(0,\infty)\backslash\cup_{k\in\mathcal{K}}\Delta_k},$$
where the collection $\{\Delta_k\}_{k\in\mathcal{K}}$ is given by Proposition \ref{Braverman-Mekler lemma}.

By Lemma \ref{general properties of fi}\eqref{second property fi}, $\varphi(z)\leq\varphi(x)=0.$ By Theorem \ref{expectation lemma}, $P(z|\mathcal{A})\in Q'(z)$ for every $\mathcal{A}\in\mathfrak{A}.$ It follows now from Corollary \ref{tech corollary} that $v\in Q'(z).$ Since $u\leq v,$ it follows from Lemma \ref{normality of qprime} that $u\in Q'(z).$ Theorem \ref{expectation lemma} implies that $z\in Q'(x).$ By Lemma \ref{transitivity lemma}, $u\in Q'(x).$ However,
$$\dist(y,Q'(x))\leq||y-u||_E\leq||y-P(y|\mathcal{A}_{\varepsilon})||_{L_1\cap L_{\infty}}\leq\varepsilon(1+y(1)).$$
Since $\varepsilon$ is arbitrary, $y\in Q'(x).$

Let now $y\in\Omega_+(x)$ be arbitrary. By \cite[Lemma II.2.1 and
Theorem II.2.1]{KPS}, for every fixed $\varepsilon>0,$  there
exist $y_1\in E,$ $y_2\in E$, $y=y_1+y_2$ and measure-preserving
transformation $\gamma$ such that $0\leq y_1\leq y^*\circ\gamma$
and $||y_2||_E\leq\varepsilon.$ Since we already proved that
$y^*\in Q'(x),$ the assertion follows immediately.
\end{proof}

\begin{thm}\label{qplus nonintegrable} Let $E=E(0,\infty)$ be a fully symmetric space on semi-axis. Suppose that $E\not\subseteq L_1$ and $x\in E.$ If $\varphi(x)=0,$ then the set $\Omega_+(x)$ is the norm-closed convex hull of its extreme points. If, in addition, the norm on $E$ is a Fatou norm, then converse assertion also holds.
\end{thm}

\begin{proof} The assertion follows immediately from Theorem \ref{qprime nonintegrable}.

By Proposition \ref{nonlinearity of fi}, there exist $0\leq y_1,z_1\in E,$ such that $x=y_1+z_1$ and $\varphi_{cut}(x)=\varphi_{cut}(y_1)=\varphi_{cut}(z_1).$By assumption, $y_1,z_1\in Q_+(x).$ By Proposition \ref{linearity of fi}, $\varphi_{cut}(x)=\varphi_{cut}(y_1)+\varphi_{cut}(z_1).$ Consequently, $\varphi_{cut}(x)=0.$ By Proposition \ref{nonlinearity of fi}, there exist $0\leq y_2,z_2\in E,$ such that $x=y_2+z_2$ and $\varphi(x)=\varphi(y_2)=\varphi(z_2).$ By Proposition \ref{linearity of fi infinite}, $\varphi(x)=\varphi(y_1)+\varphi(z_1).$ Consequently, $\varphi(x)=0.$
\end{proof}

\section{Appendix}

\subsection{An application to the case of orbits $\Omega(x)$}

The following consequence of Theorem \ref{qplus finite} is
essentially due to Braverman and Mekler \cite{BM}.

\begin{cor}\label{result of BM} If $\varphi(x)=0,$ then $\Omega(x)$ is the norm-closed convex hull of its extreme points.
\end{cor}
\begin{proof} Let $x=x^*$ and $y\in\Omega(x).$ Clearly, $y=u\cdot |y|,$ where $|u|=1$ a.e. and $|y|\in\Omega_+(x).$ Fix $\varepsilon>0.$ By Theorem \ref{qplus finite}, there exist $n\in N,$ scalars $\lambda_{n,i},\beta_{n,i}\in[0,1]$ and functions $x_{n,i}\sim x\chi_{[0,\beta_{n,i}]},$ such that $\sum_{i=1}^n\lambda_{n,i}=1$ and
$$||\,|y|-\sum_{i=1}^n\lambda_{n,i}x_{n,i}||_E\leq\varepsilon.$$
There exist measure-preserving transformations  $\gamma_{n,i},$
$1\leq i\leq n,$ (see \cite{Ryff3}) such that
$x_{n,i}=(x^*\chi_{[0,\beta_{n,i}]})\circ\gamma_{n,i}.$ Set
$x^1_{n,i}=u\cdot x\circ\gamma_{n,i}$ and
$x^2_{n,i}=u\cdot(x\chi_{[0,\beta_{n,i}]}-x\chi_{[\beta_{n,i},1]})\circ\gamma_{n,i},$
$1\leq i\leq n.$ It is clear that $x_{n,i}\sim x,$ $1\leq i\leq
n,$ and
$$||\,y-\frac12\sum_{i=1}^n\lambda_{n,i}x^1_{n,i}-\frac12\sum_{i=1}^n\lambda_{n,i}x^2_{n,i}||_E\leq\varepsilon.$$
\end{proof}

\subsection{Extreme points of the orbit $\Omega_+(x)$}

The following theorem is due to Ryff (see \cite{Ryffproc}).
\begin{thm}\label{Ryff Theorem} If $0\leq x\in L_1(0,1),$ then  $ y\in\extr(\Omega'(x))$ if and only if $y^*=x^*.$
\end{thm}

\begin{cor}  If $0\leq x\in L_1(0,1),$ then  $y\in\extr(\Omega_+(x))$ if and only if
$y^*=x^*\chi_{[0,\beta]}$ for some $\beta\geq0.$
\end{cor}
\begin{proof} Indeed, if $\int_0^{\beta}x^*(s)ds=\int_0^1y^*(s)ds,$ then
$y\in\Omega'(x^*\chi_{[0,\beta]})$. Therefore, if $y\in\extr(\Omega_+(x))$, then obviously $y\in\extr(\Omega'(x^*\chi_{[0,\beta]}))$ and the assertion
follows immediately from Theorem \ref{Ryff Theorem}.

If $y^*=x^*\chi_{[0,\beta]}$ and $y=\frac12(u_1+u_2)$ with $u_i\in\Omega_+(x),$ then $\int_0^tu^*_i(s)ds=\int_0^tx^*(s)ds$ for $t\in[0,\beta]$ and $\supp(u_i)=\supp(y).$ Therefore, $(u_1+u_2)^*=u_1^*+u_2^*.$ It follows now from \cite[(II.2.19)]{KPS} that $u_1=u_2.$
\end{proof}

\begin{lem}\label{extreme} If $0\leq x\in L_1+L_{\infty}$ and $y\in\extr(\Omega_+(x)),$ then $y\chi_{\{y<y^*(\infty)\}}=0.$
\end{lem}
\begin{proof} Assume, the contrary. Thus, the Lebesgue measure of the set
$A=\{y\in(0,\lambda y^*(\infty))\}$ does not vanish for some
$\lambda\in(0,1)$. Let $0\leq\varepsilon$ be such that
$(1+\varepsilon)\lambda<1.$ Set
$y_1=(1+\varepsilon)y\chi_A+y\chi_{(0,\infty)\backslash A}$ and
$y_2=(1-\varepsilon)y\chi_A+y\chi_{(0,\infty)\backslash A}.$
Clearly, $y_i^*=y^*$ and, therefore, $y_i\in\Omega_+(x),$ for
$i=1,2$. Hence, $y=\frac12(y_1+y_2)\notin\extr(\Omega_+(x)).$
\end{proof}

\begin{cor}\label{plus infinite} Let $0\leq x\in L_1+L_{\infty}$ and $y\in\extr(\Omega_+(x)).$ It then follows that
\begin{enumerate}
\item If $x^*(\infty)=0,$ then $y^*=x^*\chi_{[0,\beta]}$ for some $\beta\in[0,\infty].$
\item If $x^*(\infty)>0,$ then either $y^*=x^*\chi_{[0,\beta]}$ for some $\beta\in[0,\infty)$ or
$y^*=x^*$ and $y\chi_{\{y<y^*(\infty)\}}=0.$
\end{enumerate}
Conversely, functions as above belong to the set
$\extr(\Omega_+(x))$.
\end{cor}
\begin{proof} If $y$ belongs to $\extr(\Omega_+(x)),$ then so does $y^*$ (see \cite{Ryffproc} and \cite{CKS}).
Fix $t_1>0$ and find $t_2\leq t_1$ such that $\int_0^{t_2}x^*(s)ds=\int_0^{t_1}y^*(s)ds.$ Clearly, $y^*\chi_{[0,t_1]}\prec x^*\chi_{[0,t_2]}$ and $y^*\chi_{[t_1,\infty)}\prec\prec x^*\chi_{[t_2,\infty)}.$ If $y^*\chi_{[0,t_1]}=\frac12(u_1+u_2)$ with $u_1,u_2\in\Omega'(x^*\chi_{[0,t_2]}),$ then set $y_i=u_i\chi_{[0,t_1]}+y^*\chi_{[t_1,\infty)}.$ We claim $y_i\prec\prec x.$ Indeed, if $e\in(0,\infty)$ and $m(e)<\infty,$ then $e=e_1\cup e_2$ with $e_1\subset[0,t_1]$ and $e_2\subset[t_1,\infty).$ Therefore,
$$\int_ey_i(s)ds=\int_{e_1}u_i(s)ds+\int_{e_2}y^*(s)ds\leq\int_0^{m(e_1)}u_i^*(s)ds+\int_{t_1}^{t_1+m(e_2)}y^*(s)ds\leq$$
$$\leq\int_0^{\min\{t_2,m(e_1)\}}x^*(s)ds+\int_{t_2}^{t_2+m(e_2)}x^*(s)ds\leq\int_0^{m(e)}x^*(s)ds.$$
Hence, $y_i\in\Omega_+(x)$ and $y=\frac12(y_1+y_2).$ Thus,
$y\notin\extr(\Omega_+(x)).$ Therefore,
$y^*\chi_{[0,t_1]}\in\extr(\Omega'(x^*\chi_{[0,t_2]})).$ By
Theorem \ref{Ryff Theorem}, $y^*=x^*$ on $[0,t_2]$. The assertion
follows now from Lemma \ref{extreme}.

The converse assertion is easy.
\end{proof}
\begin{cor} If $x\in L_1(0,\infty),$ then $0\leq y\in\extr(\Omega'(x))$ if and only if $y^*=x^*.$
\end{cor}
The proof is identical to that of Corollary \ref{plus infinite}.

\subsection{Marcinkiewicz spaces with trivial functional $\varphi$} It follows from the Lemma \ref{general properties of fi} and the definition of Marcinkiewicz space, that $\varphi=0$ if and only if $\varphi(\psi')=0.$ It is now easy to derive, that in case of the interval $(0,1)$ this is equivalent to the condition
$$\liminf_{t\to0}\frac{\psi(2t)}{\psi(t)}>1.$$
In case of the semi-axis, the condition
$$\liminf_{t\to\infty}\frac{\psi(2t)}{\psi(t)}>1$$
needs to be added.

\subsection{A comparison of conditions \eqref{bm condition} and \eqref{new condition} in Orlicz spaces } Let $M$ be a convex function
satisfying \eqref{convex function} and let $L_M$ be the
corresponding Orlicz space on $(0,1)$. The following proposition shows that $L_M$ always satisfies condition \eqref{new condition}.

\begin{prop}\label{new condition in Orlicz} We have $\varphi(x)=0$ for every $x\in L_M.$
\end{prop}
\begin{proof}
Using the description of relatively weakly compact subsets in $L_{M}$ given in
\cite{Ando} (see also ~\cite[p.~144]{RaRe1991})
we see that for every $0\leq y\in L_M$
$$n\int_0^{\frac1n}M(\frac1ny)\to0.$$

We are going to prove that $\frac1n||\sigma_nx||_{L_M}\to0$ for
every $x\in L_M.$ Assume the contrary. Let
$||\sigma_nx||_{L_M}\geq n\alpha$ for some $0\leq x\in L_M$, some $\alpha>0$ and for arbitrary large
$n\ge 1.$ By the definition of the norm $\|\cdot \|_{L_M}$, we have
$$\int_0^1M(\frac1{n\alpha}\sigma_nx)\geq 1.$$
Hence,
$$n\int_0^{\frac1n}M(\frac1ny)\geq1$$
with $y=\alpha^{-1}x\in L_M.$ A contradiction.
\end{proof}

We shall now present an example of an Orlicz space $L_M$ which fails to satisfy condition \eqref{bm condition}.

For the definition of Boyd indices $1\leq p_{E}\leq q_{E}\leq
\infty $ of a fully symmetric space $E,$ we refer the reader to
\cite[2.b.1 and p. 132]{LT2}. It is clear, that the condition
\eqref{bm condition} holds for a fully symmetric space $E$ if and only if $p_E>1$. It is
well-known (see e.g. \cite{LT2}) that Orlicz space $L_M$ is
separable if and only if $q_{_{L_M}}<\infty$.

\begin{ex} There exists a non-separable Orlicz space $L_M$ such that $p_{L_M}=1.$
\end{ex}

\begin{proof} Let $a_0=1$ and $a_{n+1}=e^{a_n}.$ Set $M(t)=t^2$ on $(0,1),$
$M(t)=e^t+M(a_{2n})-e^{a_{2n}}$ on $[a_{2n},a_{2n+1}]$ and
$M(t)=M(a_{2n-1})+e^{a_{2n-1}}(t-a_{2n-1})$ on $[a_{2n-1},a_{2n}].$
Clearly, $M'(t)=e^t$ on $[a_{2n},a_{2n+1}]$ and $M'(t)=e^{a_{2n-1}}$ on
$[a_{2n-1},a_{2n}].$ Hence, $M'(t)\leq e^t$ and $M(t)\leq e^t-1.$

If $q_{_{L_M}}<\infty,$ then (see \cite[2.b.5]{LT2}) there exists $q$ such that
$$\sup_{\lambda,t\geq1}\frac{M(\lambda t)}{M(\lambda)t^q}<\infty.$$
In particular, $M(t)\leq const\cdot t^q$ for $t\geq1.$ However,
$$M(a_{2n+1})\geq e^{a_{2n+1}}-e^{a_{2n}}=e^{a_{2n+1}}(1+o(1)).$$
Therefore, $q_{_{L_M}}=\infty$ and $L_M$ is non-separable.

If $p_{_{L_M}}<1,$ then (see \cite[2.b.5]{LT2}) there exists $p>1$ such that
$$\inf_{\lambda,t\geq1}\frac{M(\lambda t)}{M(\lambda)t^p}>0.$$
Set $\lambda=n$ and $t=\frac1na_{2n}.$ Hence, $\lambda t=a_{2n}$ and
$$M(\lambda t)=M(a_{2n-1})+e^{a_{2n-1}}(a_{2n}-a_{2n-1})=a_{2n}(1+o(1))+a_{2n}^2(1+o(1)).$$
Since $a_{2n-1}=\frac1n o(a_{2n}),$ then
$$M(\lambda)=M(a_{2n-1})+e^{a_{2n-1}}(\frac1na_{2n}-a_{2n-1})=a_{2n}(1+o(1))+\frac1na_{2n}^2(1+o(1)).$$
Therefore,
$$\frac{M(\lambda t)}{M(\lambda)t^p}=(1+o(1))\frac{a_{2n}^2}{\frac1na_{2n}^2\cdot n^p}=(1+o(1))n^{1-p}=o(1)$$
and we conclude $p_{_{L_M}}=1$.
\end{proof}

\subsection{An application to symmetric functionals}

Let $E$ be a fully symmetric space. A positive functional $f\in E^*$ is said to be symmetric (respectively, fully symmetric) if $f(y)=f(x)$ (respectively, $f(y)\leq f(x)$) for all $0\leq x,y \in E$ such that $y^*=x^*$ (respectively, $y\prec\prec x$).
We refer to \cite{DPSS1998,CaSu} and references therein for the exposition of the theory of singular fully symmetric functionals and their applications. Recently, symmetric functionals which fail to be fully symmetric were constructed in \cite{KS} on some Marcinkiewicz spaces. However, for Orlicz spaces situation is different. The following proposition shows that a symmetric functional on an Orlicz space on the interval $(0,1)$ is necessary fully symmetric.

\begin{prop}\label{result on symm_funct} Any symmetric functional on $L_M$ is fully symmetric.
\end{prop}

\begin{proof} Let $\omega\in E^*$ be symmetric. It is clear, that $\omega(x^*\chi_{[0,\beta]})\leq\omega(x)$ for $x\geq0.$ Therefore, $\omega(y)\leq\omega(x)$ for $y\in\overline{\conv}\{y^*=x^*\chi_{[0,\beta]}\}.$ Since $\omega$ is continuous, we have $\omega(y)\leq\omega(x)$ for $y\in Q_+(x)$.
By Theorem \ref{qplus finite} and Proposition \ref{new condition in Orlicz}, we have $Q_+(x)=\Omega_+(x),$ and so $\omega$ is a fully symmetric functional on $L_M.$
\end{proof}
\begin{cor} Any singular symmetric functional on $L_M$ vanishes.
\end{cor}
\begin{proof} Indeed, there are no fully symmetric singular functionals on $L_M$ (see \cite[Theorem 3.1]{DPSS1998}).
\end{proof}

We also formulate the following hypothesis: If $E$ is a fully symmetric space, then functional $\varphi$ vanishes if and only if there are no singular symmetric functionals on $E.$

\bigskip
\leftline{F. Sukochev}\leftline{School of Mathematics and Statistics}\leftline{University of New South Wales, Kensington, NSW 2052, Australia}\leftline{Email Address:{\it f.sukochev@unsw.edu.au}}
\leftline {D. Zanin}\leftline {School of Computer Science, Engineering and Mathematics} \leftline {Flinders University,
Bedford Park, SA 5042, Australia}\leftline {Email Address: {\it
zani0005@csem.flinders.edu.au} }

\end{document}